\documentclass{amsart}
\usepackage{amsmath,amssymb,amsthm}
\usepackage{bbm}
\usepackage[left=2cm, right=2cm, top=2cm, bottom=2cm]{geometry}
\theoremstyle{definition}
\newtheorem{defn}{Definition}

\newtheorem{thm}{Theorem}
\newtheorem{theorem}[defn]{Theorem}
\newtheorem{lem}[defn]{Lemma}
\newtheorem{coro}[defn]{Corollary}
\newtheorem{prop}[defn]{Proposition}
\newtheorem{rmk}[defn]{Remark}
\newtheorem{exmp}{Example}
\newtheorem*{qn*}{Question}

\numberwithin{defn}{section}

\newcommand{\N}{\mathbb{N}}

\newcommand{\C}{\mathbb{C}}
\newcommand{\F}{\mathbb{F}}
\newcommand{\M}{\mathcal{M}}
\newcommand{\I}{\mathcal{I}}

\let\oldproofname=\proofname
\renewcommand{\proofname}{\rm\bf{\oldproofname}}
\renewcommand{\leq}{\leqslant}
\renewcommand{\geq}{\geqslant}

\makeatletter
\@namedef{subjclassname@2020}{\textup{2020} Mathematics Subject Classification}
\makeatother

\begin{document}

\pagenumbering{arabic}

\subjclass[2020]{Primary 11T55; Secondary 11N56}
\keywords{Hal\'{a}sz's Theorem, Mean values of multiplicative functions, Arithmetic of polynomials over finite fields}

\title{A note on Hal\'{a}sz's Theorem in $\mathbb{F}_q[t]$}
\author[A. Afshar]{Ardavan Afshar}
\address{Department of Mathematics \\ KTH \\
Lindstedtsv\"{a}gen 25 \\ Stockholm \\ Sweden}
\email{ardavana@kth.se}

\maketitle

\begin{abstract}
In the setting of the integers, Granville, Harper and Soundararajan showed that the upper bound in Hal\'{a}sz's Theorem can be improved for smoothly supported functions. We derive the analogous result for Hal\'{a}sz's Theorem in $\mathbb{F}_q[t]$, and then consider the converse question of when the general upper bound in this version of Hal\'{a}sz's Theorem is actually attained. 
\end{abstract}

\section{Introduction}

\subsection{Hal\'{a}sz's Theorem for the integers}

Let $f: \N \to \C$ be a multiplicative function such that $f(1)=1$, and such that its associated Dirichlet Series and Euler Product (respectively)
$$ \mathcal{F}(s) := \sum_{n \geq 1} \frac{f(n)}{n^s} = \prod_{p \text{ prime }} \sum_{k \geq 0} \frac{f(p^k)}{p^{ks}} $$
are defined and absolutely convergent for $s \in \C$ with $\Re(s) > 1$. Then define $\Lambda_f(n)$, the von Mangoldt function associated to $f$, by
$$ - \frac{\mathcal{F}'}{\mathcal{F}} (s) =: \sum_{n \geq 1} \frac{\Lambda_{f}(n)}{n^s}$$
and consider the set of such functions $f$ such that, for some $\kappa > 0$, we have $ |\Lambda_{f}(n)| \leq  \kappa \Lambda(n)$ for all $n \geq 1$ (where $\Lambda$ is the usual von Mangoldt function), which we denote $\mathcal{C}(\kappa)$. In \cite{LongPaper}, Granville, Harper and Soundararajan generalise Hal\'{a}sz's Theorem to this class of functions:

\begin{theorem}[Hal\'{a}sz's Theorem] \label{Halasz}
Let $\kappa > 0$, $x$ large and $f \in \mathcal{C}(\kappa)$, and define $M = M(x)$ by
$$ e^{-M} (\log x)^\kappa := \max_{|t| \leq (\log x)^\kappa} \left \lvert  \frac{\mathcal{F}\left(1+1/\log x + it\right)}{1 + 1/\log x + it} \right \rvert. $$
Then we have that
$$ S(x) := \frac{1}{x} \left|\sum_{n \leq x} f(n)\right| \ll_\kappa (1+M)e^{-M}(\log x)^{\kappa-1} + \frac{(\log \log x)^\kappa}{\log x}.$$
\end{theorem}

\begin{rmk}
Hal\'{a}sz's Theorem gives us a very general tool for understanding multiplicative functions, and provides another way to recover results associated to particular cases. Note, for example, that the non-vanishing of $\zeta(s)$ on $\Re(s)=1$ implies, for $f = \mu$, the M\"{o}bius function, that $e^{-M} \asymp \frac{1}{\log x}$ and so by Hal\'{a}sz's Theorem we have
$$ \frac{1}{x} \left| \sum_{n \leq x} \mu(n) \right| \ll \frac{\log \log x}{\log x}$$
which is equivalent to the Prime Number Theorem (albeit with a weak error term).
\end{rmk}

\begin{rmk} \label{SmoothIntegers}
In the case of $f \in \mathcal{C}(1)$ the inequality in Theorem \ref{Halasz} becomes
\begin{equation} \label{eqn:Halasz}
S(x) \ll (1+M)e^{-M} + \frac{(\log \log x)}{\log x} .
\end{equation}
Now, for simplicity, consider the multiplicative functions $f$ with $f(1)=1$ and $|f(n)| \leq 1$ for all $n$, which form a superset of $\mathcal{C}(1)$. For this set, the same authors show that if $f$ is supported only on primes of size $p \leq x^{1 - \delta}$ for some $\delta > 0$, then we can improve the upper bound in equation \eqref{eqn:Halasz} to get
\begin{equation} \label{eqn: smoothbound}
S(x) \ll_\delta e^{-M} + \frac{(\log \log x)}{\log x}.
\end{equation}
This observation is presented in Remark 3.2 of \cite{ShortPaper}, albeit with a different set of notation associated to this setting.
\end{rmk}

\begin{rmk}
The upper bound in \eqref{eqn: smoothbound} does not hold for general multiplicative functions. It has been shown by a variety of authors (see \cite{M}, \cite{M-V} and \cite{Decay}),  following the idea of Montgomery in \cite{M}, that there exists a multiplicative function $f$ with $f(1)=1$ and $|f(n)| \leq 1$ for all $n$, such that
\begin{equation}
S(x) \gg (1+M)e^{-M} + \frac{(\log \log x)}{\log x}.
\end{equation}
\end{rmk}

\subsection{Hal\'{a}sz's Theorem in $\mathbb{F}_q[t]$}

We work in the setting of polynomials over a finite field, and set up the quantities analogous to those in the setting of the integers by following the notation in \cite{FunctionFieldsPaper}. Let $\F_q$ be a finite field of order $q$, $\M = \{F \in \F_q[t] \text{ monic}\}$ and $\I= \{P \in \M : P \text{ is irreducible}\}$. We define $f: \M \to \C$ to be a multiplicative function such that $f(1)=1$, and such that its associated Power Series and Euler Product (respectively)
$$ \mathcal{F}(z) := \sum_{F \in \M} f(F)z^{\deg F}  = \prod_{P \in \I}\sum_{k \geq 0} f(P^k)z^{k \deg P} $$
are defined and absolutely convergent for $z \in \C$ with $|z| < \frac{1}{q}$. By taking the logarithmic derivative of the latter, and multiplying by $z$, we acquire a new power series through which we can define $\Lambda_f(F)$ 
(the von Mangoldt function associated to $f$):
$$ \frac{z\mathcal{F}'}{\mathcal{F}}(z) =: \sum_{F \in \M} \Lambda_f(F)z^{\deg F} .$$
Then, we let $\M_n = \{F \in \M : \deg F = n\}$, and define
$$ \sigma(n) = \sigma(n; f) : = \frac{1}{q^n} \sum_{F \in \M_n} f(F) $$
to be the mean value of $f$ over polynomials of degree $n$ and
$$ \chi(n) = \chi(n; f) : = \frac{1}{q^n} \sum_{F \in \M_n} \Lambda_f(F) $$
to be the corresponding weighted average over prime powers. \\

As in the setting of the integers, we consider the set $\mathcal{C}(\kappa)$ of such $f$ such that, for some $\kappa > 0$, we have that $ |\Lambda_{f}(F)| \leq  \kappa \Lambda(F)$ for all $F \in \M$, where
$$ \Lambda(f) = \begin{cases} \deg P &\text{ if } F = P^k \\ 0 & \text{else}\end{cases} .$$ 
In particular, given the prime polynomial theorem in the form $\sum_{F \in \M_n} \Lambda_{f}(F) = q^n$, for $f \in \mathcal{C}(\kappa)$ we have that
$$ |\chi(n)| \leq \frac{1}{q^n} \sum_{F \in \M_n} |\Lambda_f(F)| \leq \kappa $$  
and so we consider the more general set $\tilde{\mathcal{C}}(\kappa)$ of $f$ with $ |\chi(j)| = |\chi(j; f)| \leq \kappa $ for all $j \geq 1$. \\

Finally, we define $f^\bot = f^{\bot, n}$ by setting 
$$\Lambda_{f^\bot}(F) = \begin{cases} \Lambda_f(F) &\text{ if } \deg F < n \\ 0 &\text{ else} \end{cases} $$
and then we set $ \mathcal{F^\bot}(z) := \sum_{F \in \M} f^{\bot}(F)z^{\deg F} $, $ \sigma^\bot(j)  : = \frac{1}{q^j} \sum_{F \in \M_j} f^\bot(F) $ and $ \chi^\bot(j): = \frac{1}{q^n} \sum_{F \in \M_j} \Lambda_{f^\bot}(F) $. We observe that $\chi^\bot(j) = \chi(j)$ if $ j < n $ and $\chi^\bot(j) = 0$ otherwise, and from equation (1.8) of \cite{FunctionFieldsPaper} we have
\begin{equation} \label{eqn:convolution}
j\sigma(j) =  \sum_{k=1}^j \chi(k) \sigma(j-k)
\end{equation}
from which we conclude that $\sigma^\bot(j) = \sigma(j)$ if  $ j < n $ and $\sigma^\bot(n) = \sigma(n) - \frac{\chi(n)}{n}$. \\

We also note that, from their definitions and our observations above, we have
\begin{equation} \label{eqn:F}
 \mathcal{F}\left(\frac{z}{q}\right) = \sum_{j \geq 0} \sigma(j) z^j = \exp\left(\sum_{j\geq 1} \frac{\chi(j)}{j} z^j\right)
\end{equation}
and
\begin{equation} \label{eqn:Fbot}
 \mathcal{F^\bot}\left(\frac{z}{q}\right) = \sum_{j \geq 0} \sigma^\bot(j) z^j = \exp\left(\sum_{j\geq 1} \frac{\chi^\bot(j)}{j} z^j\right) = \exp\left(\sum_{j=1}^{n-1} \frac{\chi(j)}{j} z^j\right).
\end{equation}

With these definitions in place, we are able to formulate the analogue of Hal\'{a}sz's Theorem in {$\F_q[t]$}, which Granville, Harper and Soundararajan prove in \cite{FunctionFieldsPaper}: 
\begin{theorem}[Hal\'{a}sz's Theorem in {$\F_q[t]$}] \label{HalaszFF}
Let $\kappa > 0$, $n \geq 1$ and $f \in \tilde{\mathcal{C}}(\kappa)$, and define $M = M(n)$ by
$$ e^{-M} (2n)^\kappa := \max_{|z| = \frac{1}{q}} \lvert  \mathcal{F}^\bot(z) \rvert. $$
Then we have that
\begin{equation} \label{eqn:HalaszFF}
|\sigma(n)|  \leq 2\kappa(\kappa + 1+M)e^{-M}(2n)^{\kappa-1}.
\end{equation}
\end{theorem}

\begin{rmk}
In Theorem \ref{Halasz}, we define $M$ in terms of the maximum value of the Dirichlet Series $\mathcal{F}(s)$ on the line segment $\{\Re(s) = 1 + 1/\log x + it : |t| \leq (\log x)^\kappa \}$. The restriction of $t$ up to height $(\log x)^\kappa$ comes from using a truncated version of Pellet's formula in the proof of Theorem \ref{Halasz}, and taking the real part of $s$ to be $1 + 1/\log x$ is the ensure the convergence of $\mathcal{F}(s)$. The analogous proof of Theorem \ref{HalaszFF} uses Cauchy's theorem, in which we integrate over the whole circle, and to ensure convergence, instead truncates the Power Series $\mathcal{F}(z)$ at height $n$ (which is equivalent to its analogue, up to a multiplicative constant). This is why, in Theorem \ref{HalaszFF}, we define $M$ in terms of the maximum value of the Power Series $\mathcal{F}^\bot(z)$ on the circle $|z| = 1/q$.
\end{rmk}

We consider the case analogous to that discussed in Remark \ref{SmoothIntegers}, and show that the upper bound in Hal\'{a}sz's Theorem can be improved when $f$ is smoothly supported.
\begin{thm} \label{HalaszFFSmooth}
Let $\kappa > 0$, $n \geq 1$ and $f \in \tilde{\mathcal{C}}(\kappa)$, and define $M = M(n)$ as in Theorem \ref{HalaszFF}. Suppose in addition that, for some small $\delta > 0$,  $f$ is supported only on irreducibles $P$ of degree at most $(1-\delta)(n-1)$.
Then we get that
\begin{align*}
|\sigma(n)| &\leq 2\kappa^2 e^{-M} (2n)^{\kappa - 1} \left(1 - \log \left(1 - e^{-\frac{\delta}{2\sqrt{1-\delta}}}\right) + \frac{1}{\kappa}\left(1-\frac{e^{M/\kappa}}{2n}\right)^{\delta(n-1)} \right) + \left(\frac{q}{q-1}\right)^2\frac{\kappa n^{\kappa-2}}{q^{(1-\delta)(n-1)/2}}\\
&\ll_\delta \kappa(\kappa + 1)e^{-M} (2n)^{\kappa - 1}. 
\end{align*}
\end{thm}

Conversely, we derive a criterion for when the upper bound in equation \eqref{eqn:HalaszFF} is asymptotically attained:

\begin{thm} \label{converse}
Let $\kappa > 0$, $n \geq 1$ and $f \in \tilde{\mathcal{C}}(\kappa)$, and define $M = M(n)$ as in Theorem \ref{HalaszFF}. Suppose that $\kappa + 1 = o(M)$, then
$ |\sigma(n)| \gg \kappa(\kappa + 1 + M)e^{-M} (2n)^{\kappa - 1} $
if, and only if, for all $\delta \gg 1$ we have
\begin{equation*} \label{eqn:fund}
\left|\sum_{(1-\delta)(n-1) < j \leq n} \chi(j) \sigma(n-j)\right| \gg \kappa M e^{-M}(2n)^{\kappa}.
\end{equation*}
\end{thm}
\begin{rmk}
The assumption that $\kappa + 1 = o(M)$ in Theorem \ref{converse} is precisely the case in which the upper bound in Theorem \ref{HalaszFFSmooth} is actually asymptotically smaller than the upper bound in Theorem \ref{HalaszFF}.
\end{rmk}

Finally, inspired by the idea in the setting of the integers (see \cite{M}), in section \ref{sharpexample} we compute an example for which the criterion in Theorem \ref{converse} holds, in the case of $\kappa = 1$. \\

\section{Proofs of Theorems \ref{HalaszFFSmooth} and \ref{converse}}

Let $\kappa > 0$, and let $f \in \tilde{\mathcal{C}}(\kappa)$. From equation (3.3) of \cite{FunctionFieldsPaper} we have that 
\begin{equation}
\sigma(n) - \frac{\chi(n)}{n}= \frac{1}{nq^n} \int_0^1 \frac{1}{2\pi i}\int_{|z| = \frac{1}{q\sqrt{t}}} \left(\sum_{j=1}^{n-1} \chi(j)(qz)^j \right) \left(\sum_{j=1}^{n-1} \chi(j)(qtz)^j \right)\mathcal{F^\bot}(tz)\frac{dz}{z^{n+1}} \frac{dt}{t}.
\end{equation}
We define a new quantity for $n > m \geq 1$
\begin{equation} \label{eqn:sigma_m}
\sigma_m(n) := \frac{1}{nq^n} \int_0^1 \frac{1}{2\pi i}\int_{|z| = \frac{1}{q\sqrt{t}}} \left(\sum_{j=1}^{m} \chi(j)(qz)^j \right) \left(\sum_{j=1}^{n-1} \chi(j)(qtz)^j \right)\mathcal{F^\bot}(tz)\frac{dz}{z^{n+1}} \frac{dt}{t}
\end{equation}
so that $\sigma_{n-1}(n) = \sigma(n) - \frac{\chi(n)}{n}$, and bound it following the strategy in \cite{FunctionFieldsPaper}.

\begin{prop}
Let $\kappa > 0$, $n \geq 1$, $f \in \tilde{\mathcal{C}}(\kappa)$ and $M = M(n)$ as in Theorem \ref{HalaszFF}. Then for $m < n-1$ we have
$$ |\sigma_m(n)| \leq 2\kappa^2 e^{-M} (2n)^{\kappa - 1} \left(1 - \log \left(1 - e^{-\frac{(n-1)-m}{2\sqrt{m(n-1)}}}\right) + \frac{1}{\kappa}\left(1-\frac{e^{M/\kappa}}{2n}\right)^{(n-1)-m} \right). $$
\end{prop}
\begin{proof}
First we use Cauchy-Schwarz on the inner integral in equation \eqref{eqn:sigma_m}
\begin{align*}
&\left|\frac{1}{2\pi i}\int_{|z| = \frac{1}{q\sqrt{t}}} \left(\sum_{j=1}^{m} \chi(j)(qz)^j \right) \ \left(\sum_{j=1}^{n-1} \chi(j)(qtz)^j \right)\mathcal{F^\bot}(tz)\frac{dz}{z^{n+1}} \right| \\
&\leq (q\sqrt{t})^n \left(\max_{|z| = \frac{1}{q\sqrt{t}}} |\mathcal{F^\bot}(tz)|\right)\sqrt{I_m\left(1, \frac{1}{q\sqrt{t}}\right)I_{n-1}\left(t, \frac{1}{q\sqrt{t}}\right)}
\end{align*}
where, for $ s \geq 0$, we have
$$I_a(s, R) := \frac{1}{2\pi} \int_{|z| = R} \left|\sum_{j=1}^{a} \chi(j)(qsz)^j \right|^2 \frac{|dz|}{|z|} =  \sum_{j=1}^{a} |\chi(j)|^2 \left(qsR\right)^{2j} \leq \kappa^2 \sum_{j=1}^{a} \left(qsR\right)^{2j} $$ 
by Parseval's identity. Using this, we bound the inner integral by the quantity
\begin{align*}
&\kappa^2 (q \sqrt{t})^n \left(\max_{|z| = \frac{1}{q\sqrt{t}}} |\mathcal{F^\bot}(tz)|\right) \left(\sum_{j=1}^{m} t^{-j} \right)^{\frac{1}{2}} \left(\sum_{j=1}^{n-1}t^j \right)^{\frac{1}{2}}  = \kappa^2 q^n\left(\max_{|z| = \frac{\sqrt{t}}{q}} |\mathcal{F^\bot}(z)|\right) t^{\frac{n-m+1}{2}}\left(\frac{\sqrt{(1-t^{m})(1-t^{n-1})}}{1-t}\right) 
\end{align*}
and then recall the bound from equation (3.6) of \cite{FunctionFieldsPaper}, which for $t \in (0,1)$, states that 
$$ \max_{|z| = \frac{\sqrt{t}}{q}} |\mathcal{F^\bot}(z)| \leq \min(e^{-M} (2n)^\kappa, (1-\sqrt{t})^{-\kappa})$$
where $ e^{-M} (2n)^\kappa :=  \max_{|z| = \frac{1}{q}} |\mathcal{F^\bot}(z)| $. \\

Putting this all back into the full integral we get
\begin{align*}
|\sigma_m(n)| &\leq \frac{\kappa^2}{n} \int_0^1  \min(e^{-M} (2n)^\kappa, (1-\sqrt{t})^{-\kappa}) \ t^\frac{(n-1)-m}{2}\left(\frac{\sqrt{(1-t^{m})(1-t^{n-1})}}{1-t}\right) dt
\end{align*}
and after the substitution $t=(1-u)^2$ we have
\begin{align*}
|\sigma_m(n)| &\leq \frac{\kappa^2}{n} \int_0^1 \min(e^{-M} (2n)^\kappa, u^{-\kappa}) \ (1-u)^{(n-1)-m} \min\left(\sqrt{m(n-1)}, \frac{1}{u(2-u)}\right) 2(1-u) \ du \\
&\leq \frac{\kappa^2}{n} \int_0^1 \min(e^{-M} (2n)^\kappa, u^{-\kappa}) \ (1-u)^{(n-1)-m} \min\left(2\sqrt{m(n-1)}, \frac{1}{u}\right) du.
\end{align*}
Now, if $e^{M/\kappa} \sqrt{m(n-1)} \geq n$, we get
\begin{align*}
&|\sigma_m(n)| \leq \frac{\kappa^2}{n} \left( \int_0^\frac{1}{2\sqrt{m(n-1)}} 2\sqrt{m(n-1)}e^{-M} (2n)^\kappa  du + \int_{\frac{1}{2\sqrt{m(n-1)}}}^{\frac{e^{M/\kappa}}{2n}} e^{-M} (2n)^\kappa (1-u)^{(n-1)-m} \frac{du}{u} \right. \\
&+ \left. \int_{\frac{e^{M/\kappa}}{2n}}^1 (1-u)^{(n-1)-m} \frac{du}{u^{\kappa+1}} \right)
\end{align*}
and otherwise
\begin{align*}
&|\sigma_m(n)| \leq \frac{\kappa^2}{n} \left( \int_0^\frac{1}{2\sqrt{m(n-1)}} 2\sqrt{m(n-1)}e^{-M} (2n)^\kappa du + \int_{\frac{e^{M/\kappa}}{2n}}^1 (1-u)^{(n-1)-m} \frac{du}{u^{\kappa+1}} \right).
\end{align*}
We can combine these two cases as follows
\begin{align*}
&|\sigma_m(n)|\leq \frac{\kappa^2}{n} \left(e^{-M} (2n)^\kappa + e^{-M} (2n)^\kappa \sum_{j=1}^{\lceil e^{M/\kappa}\frac{\sqrt{m(n-1)}}{n} - 1\rceil} \int_\frac{j}{2\sqrt{m(n-1)}}^\frac{j+1}{2\sqrt{m(n-1)}} (1-u)^{(n-1)-m} \frac{du}{u} \right. \\
&\left.+  \left(1-\frac{e^{M/\kappa}}{2n}\right)^{(n-1)-m} \int_{\frac{e^{M/\kappa}}{2n}}^1 \frac{du}{u^{\kappa+1}}  \right).
\end{align*}

When $m < n-1$ we get
\begin{align*}
|\sigma_m(n)|
&\leq 2\kappa^2 e^{-M} (2n)^{\kappa - 1} \left(1 + \sum_{j \geq 1} \left(1-\frac{j}{2\sqrt{m(n-1)}}\right)^{(n-1)-m} \log\left(1 + \frac{1}{j}\right)  + \frac{1}{\kappa}\left(1-\frac{e^{M/\kappa}}{2n}\right)^{(n-1)-m} \right) \\
&\leq 2\kappa^2 e^{-M} (2n)^{\kappa - 1} \left(1 + \sum_{j \geq 1}  \frac{e^{-\frac{j((n-1)-m)}{2\sqrt{m(n-1)}}}}{j} + \frac{1}{\kappa}\left(1-\frac{e^{M/\kappa}}{2n}\right)^{(n-1)-m} \right) \\
&\leq 2\kappa^2 e^{-M} (2n)^{\kappa - 1} \left(1 - \log \left(1 - e^{-\frac{(n-1)-m}{2\sqrt{m(n-1)}}}\right) + \frac{1}{\kappa}\left(1-\frac{e^{M/\kappa}}{2n}\right)^{(n-1)-m} \right).
\end{align*}
\end{proof}
\begin{coro} \label{delta}
Let $\kappa > 0$, $n \geq 1$, $f \in \tilde{\mathcal{C}}(\kappa)$ and $M = M(n)$ as in Theorem \ref{HalaszFF}. Then for $\delta > 0$ and $m \leq (1-\delta)(n-1)$ we have that
\begin{align*}
|\sigma_m(n)| &\leq 2\kappa^2 e^{-M} (2n)^{\kappa - 1} \left(1 - \log \left(1 - e^{-\frac{\delta}{2\sqrt{1-\delta}}}\right) + \frac{1}{\kappa}\left(1-\frac{e^{M/\kappa}}{2n}\right)^{\delta(n-1)} \right).
\end{align*}
\end{coro}

Then we relate our quantity $\sigma_m(n)$ to $\sigma(n)$ with the following observation
\begin{lem} \label{sigmaRelation} Let $n > m \geq 1$. Then we have that
$$ \sigma(n) = \sigma_m(n) + \frac{1}{n} \sum_{j=m+1}^n \chi(j) \sigma(n-j) .$$
\end{lem}
\begin{proof}
From the definition of $\sigma_m(n)$ in equation \eqref{eqn:sigma_m}, and our observation in equation \eqref{eqn:Fbot} we have
\begin{align*}
\sigma_m(n) &= \frac{1}{nq^n} \int_0^1 \frac{1}{2\pi i}\int_{|z| = \frac{1}{q\sqrt{t}}} \left(\sum_{j=1}^{m} \chi(j)(qz)^j \right) \left(\sum_{k=1}^{n-1} \chi(k)(qtz)^k \right)\left(\sum_{l \geq 0} \sigma^\bot(l)  (qtz)^l \right)\frac{dz}{z^{n+1}}   \frac{dt}{t} \\
&= \frac{1}{nq^n} \int_0^1 \frac{1}{2\pi i}\int_{|z| = \frac{1}{q\sqrt{t}}} \left(\sum_{j=1}^{m} \chi(j)(qz)^j \right) \left(\sum_{k=1}^{n-1} \chi(k)(qtz)^k \right)\left(\sum_{l = 0}^{n-1} \sigma(l)  (qtz)^l \right)\frac{dz}{z^{n+1}}   \frac{dt}{t} \\
&= \frac{1}{nq^n} \int_0^1 \sum_{\substack{j + k + l = n \\ j \leq m}}  \chi(j)  \chi(k) \sigma(l) q^{j+k+l} t^{k+l} \frac{dt}{t} \\
&= \frac{1}{n} \sum_{\substack{j + k + l = n \\ j \leq m}}  \frac{\chi(j)  \chi(k) \sigma(l)}{k+l} = \frac{1}{n} \sum_{j=1}^m \chi(j) \frac{1}{n-j} \sum_{k=1}^{n-j} \chi(k) \sigma(n-j-k) = \frac{1}{n} \sum_{j=1}^m \chi(j) \sigma(n-j)
\end{align*}
where in the second equality we use the fact that $\sigma^\bot(j) = \sigma(j)$ for $ j < n $, and in the final equality we use equation \eqref{eqn:convolution}. Finally, using equation \eqref{eqn:convolution} once more we get that
$$ \sigma(n) = \frac{1}{n} \sum_{j=1}^n \chi(j) \sigma(n-j) = \sigma_m(n) + \frac{1}{n} \sum_{j=m+1}^n \chi(j) \sigma(n-j).$$
\end{proof}

This bring us to our proof of Theorem \ref{HalaszFFSmooth},

\begin{proof}[Proof of Theorem \ref{HalaszFFSmooth}]
Let $m = \lfloor (1-\delta)(n-1) \rfloor$. Since $f$ is supported only on irreducibles $P$ of degree at most $m$, we have for $j > m$, that
$$ |\chi(j)| = \frac{1}{q^j}\left|\sum_{F \in \M_j} \Lambda_f(F)\right| \leq \frac{1}{q^j} \sum_{\substack{d | j \\ d \leq m}} \left| \sum_{P \in \I_d} \Lambda_f(P^{j/d})\right| \leq \frac{1}{q^j} \sum_{d \leq \min(j/2, m)} q^d |\chi(d)| \leq \kappa \min(q^{-j/2}, q^{m-j}) \frac{q}{q-1} .$$
Moreover, using equation \eqref{eqn:convolution}, and our assumption that $|\chi(j)| \leq \kappa$ for all $j$, we can deduce inductively (given the base case $\sigma(0)=1$)  the trivial bound $|\sigma(j)| \leq (j+1)^{\kappa-1}$ for all $j$. Now we can bound the following sum thus
$$ \left| \sum_{j=m+1}^n \chi(j) \sigma(n-j) \right| \leq \kappa n^{\kappa-1} \frac{q}{q-1 } \sum_{j=m+1}^n q^{-j/2} \leq \left(\frac{q}{q-1}\right)^2\frac{\kappa n^{\kappa-1}}{q^{(m+1)/2}} \leq \left(\frac{q}{q-1}\right)^2\frac{\kappa n^{\kappa-1}}{q^{(1-\delta)(n-1)/2}}.$$
We combine this with Lemma \ref{sigmaRelation} and Corollary \ref{delta} to get that
\begin{equation} \label{eqn:main}
|\sigma(n)| \leq 2\kappa^2 e^{-M} (2n)^{\kappa - 1} \left(1 - \log \left(1 - e^{-\frac{\delta}{2\sqrt{1-\delta}}}\right) + \frac{1}{\kappa}\left(1-\frac{e^{M/\kappa}}{2n}\right)^{\delta(n-1)} \right) + \left(\frac{q}{q-1}\right)^2\frac{\kappa n^{\kappa-2}}{q^{(1-\delta)(n-1)/2}}.
\end{equation}
Finally, by the maximum modulus principle, $e^{-M}(2n)^{\kappa} = \max_{|z| = \frac{1}{q}}|\mathcal{F^\bot}(z)| \geq |\mathcal{F^\bot}(0)|=1$, which means that the second term in \eqref{eqn:main} is much smaller (asymptotically in $n$) than the first. So, we can conclude that
$$|\sigma(n)| \ll_\delta \kappa(\kappa + 1)e^{-M} (2n)^{\kappa - 1}.$$
\end{proof}

and our proof of Theorem \ref{converse}

\begin{proof}[Proof of Theorem \ref{converse}]
Let $\delta \gg 1$ and  $m = \lfloor (1-\delta)(n-1) \rfloor$, so that, by Lemma \ref{sigmaRelation} we have that
$$ \sigma(n) = \sigma_m(n) + \frac{1}{n} \sum_{j=m+1}^n \chi(j) \sigma(n-j) = \sigma_m(n) + \frac{1}{n}\sum_{(1-\delta)(n-1) < j \leq n} \chi(j) \sigma(n-j)$$
and from Corollary \ref{delta} we know that $\sigma_m(n) \ll \kappa(\kappa + 1) e^{-M} (2n)^{\kappa - 1}$. Therefore,  if $\kappa + 1 = o(M)$, we have that $ |\sigma(n)| \gg \kappa(\kappa+1+M)e^{-M} (2n)^{\kappa - 1} $
if, and only if,
\begin{equation} \label{eqn:converse}
\left|\sum_{(1-\delta)(n-1) < j \leq n} \chi(j) \sigma(n-j)\right| \gg \kappa Me^{-M}(2n)^{\kappa}.
\end{equation}
\end{proof}

\section{A sharp example} \label{sharpexample}

We conclude with an example for which the criterion in Theorem \ref{converse} holds, and thus which attains the upper bound in Hal\'{a}sz's Theorem. For simplicity, we take $\kappa = 1$ throughout this example.

\begin{rmk} \label{converseRemark}
First we observe that, if $0 < \delta < \frac{1}{2} - \frac{1}{2n}$, then the values taken by $\sigma(n-j)$ in equation \eqref{eqn:converse} are independent of those taken by $\chi(j)$. So, we may for example take $\chi(j) = e^{i (\theta - \phi_{n-j})}$ where $\sigma(j) =: |\sigma(j)| e^{i \phi_j}$ (for some $\theta \in [0, 2\pi)$) for $j > (1-\delta)(n-1)$. In this case Theorem \ref{converse} says that if $1=o(M)$ then
$$|\sigma(n)| \gg (1+M)e^{-M} \Leftrightarrow \sum_{j < 1 + \delta (n-1)} |\sigma(j)| \gg Me^{-M}n.  $$
\end{rmk}

We use this observation to construct the following example

\begin{exmp} Let $n \geq 2$, $0 < \delta < \frac{1}{2} - \frac{1}{2n}$ and $M = M(n)$ as in Theorem \ref{HalaszFF}. Let
$$ \chi(j) = \begin{cases} i &\text{ if } 1 \leq j < 1 + \delta (n-1) \\ 0 &\text{ if } 1 + \delta (n-1) \leq j \leq (1-\delta)(n-1) \\ e^{-i \phi_{n-j}} &\text{ if } j > (1-\delta)(n-1) \end{cases} $$
where $\sigma(j) =: |\sigma(j)| e^{i \phi_j}$. Then $1 = o(M)$ and
$$  \sum_{j < 1 + \delta (n-1)} |\sigma(j)| \gg Me^{-M}n$$
which means that, by Theorem \ref{converse}, we have
$$ |\sigma(n)| \gg (1+M)e^{-M} .$$
\end{exmp}
\begin{proof}
In this case, we have from equation \eqref{eqn:Fbot} that
\begin{align*}
\max_{|z| = \frac{1}{q}} \log|\mathcal{F^\bot}(z)| &=  \max_{|z| = \frac{1}{q}} \Re\left(\sum_{j=1}^{n-1} \frac{\chi(j)}{j}(qz)^j\right) \\
 &= \max_{\theta \in [0, 2 \pi)}\left(\sum_{1 \leq j < 1 + \delta (n-1)} \frac{-\sin(j \theta)}{j} + \sum_{(1-\delta)(n-1) < j \leq n-1} \frac{\cos(j(\theta-\phi_{n-j}))}{j} \right).
\end{align*}

Now, we know that uniformly for $\theta$ and $x$, we have $ \left| \sum_{j \leq x} \frac{\sin(j \theta)}{j}\right| \ll 1$ and moreover

\begin{align*}
\left| \sum_{(1-\delta)(n-1) < j \leq n-1} \frac{\cos(j(\theta-\phi_{n-j}))}{j} \right| &\leq \sum_{(1-\delta)(n-1) < j \leq n-1} \frac{1}{j} \ll -\log(1-\delta) \ll 1.
\end{align*}

Therefore, we have that  $ \max_{|z| = \frac{1}{q}} \log|\mathcal{F^\bot}(z)| \ll 1$, and conversely, by the maximum modulus principle $$\max_{|z| = \frac{1}{q}}|\mathcal{F^\bot}(z)| \geq |\mathcal{F^\bot}(0)| = 1.$$ This means that $ e^{-M}(2n) := \max_{|z| = \frac{1}{q}}|\mathcal{F^\bot}(z)| \asymp 1  $ and $ M = \log 2n -\max_{|z| = \frac{1}{q}} \log|\mathcal{F^\bot}(z)| \sim \log n $ so that overall we get have that $1=o(M)$ and $Me^{-M}n \asymp \log n$. \\

On the other hand, by Cauchy's Theorem, we have for  $ j < 1 + \delta (n-1)$ and $R<1$ that
\begin{align*}
\sigma(j) = \frac{1}{q^j} \frac{1}{2\pi i} \int_{|z| = \frac{R}{q}} \mathcal{F}(z)\frac{dz}{z^{j+1}}  &= \frac{1}{2\pi i} \int_{|w| = R} \mathcal{F}\left(\frac{w}{q}\right)\frac{dw}{w^{j+1}} \\
&= \frac{1}{2\pi i} \int_{|w| = R} \exp\left(\sum_{k \geq 1} \frac{\chi(k)}{k} w^k\right)\frac{dw}{w^{j+1}} \\
&= \frac{1}{2\pi i} \int_{|w| = R} \exp\left(\sum_{k=1}^j \frac{\chi(k)}{k} w^k\right)\frac{dw}{w^{j+1}} \\
&= \frac{1}{2\pi i} \int_{|w| = R} \exp\left(i \sum_{k=1}^j \frac{w^k}{k} \right)\frac{dw}{w^{j+1}} \\
&= \frac{1}{2\pi i} \int_{|w| = R} \exp\left(i \sum_{k \geq 1} \frac{w^k}{k} \right)\frac{dw}{w^{j+1}} \\
&= \frac{1}{2\pi i} \int_{|w| = R} \frac{1}{(1-w)^i}\frac{dw}{w^{j+1}} = \binom{i + j - 1}{j} \sim \frac{j^{i-1}}{\Gamma(i)}
\end{align*}
where we use equation \eqref{eqn:Fbot} in the third line. From this we conclude that
$$ \sum_{j < 1 + \delta (n-1)} |\sigma(j)| \asymp \sum_{j < 1 + \delta (n-1)} \frac{1}{j} \gg_\delta \log n \gg Me^{-M}n. $$
\end{proof}

\section*{Acknowledgements}
I would like to thank Andrew Granville for useful discussions, contributions and references, and P\"{a}r Kurlberg for his encouragement and thoughtful advice. The research leading to these results has received funding from the European Research Council (grant n$^{\text{o}}$ 670239) and the Swedish Research Council (grant n$^{\text{o}}$ 2016-03701). This work is an updated version of chapter 5 of the author's PhD thesis \cite{afshar}.

\bibliographystyle{plain}
\bibliography{example}

\end{document}